\documentclass[12pt]{amsart}
\usepackage{amssymb,stmaryrd,xypic,mathdots,graphicx}
\usepackage[all]{xy}

\theoremstyle{theorem}
\newtheorem{theorem}{Theorem}[section]
\newtheorem{lemma}[theorem]{Lemma}
\newtheorem{proposition}[theorem]{Proposition}
\newtheorem{corollary}[theorem]{Corollary}

\theoremstyle{definition}
\newtheorem{definition}[theorem]{Definition}
\newtheorem{remark}[theorem]{Remark}
\newtheorem{question}[theorem]{Question}

\title[Nonlocal loss of first homotopy]{\large N\lowercase{onlocal loss of first homotopy in polyhedral approximations of} P\lowercase{eano continua}}

\author{Jeremy Brazas}

\address{Department of Mathematics, West Chester University of Pennsylvania, West Chester, PA 19383, USA}

\email{jbrazas@wcupa.edu}

\author{Hanspeter Fischer}

\address{Department of Mathematical Sciences, Ball State University, Muncie, IN 47306, USA}

\email{hfischer@bsu.edu}

\date{August 1, 2025. \\ \mbox{\hspace{5pt} } 2020 {\em Mathematics Subject Classification.} Primary 55Q52, 55Q07; Secondary 55Q05, 54F50}


\begin{document}

\begin{abstract}
If a Peano continuum $X$ is semilocally simply connected, then it has a finite polyhedral approximation whose fundamental group is isomorphic to that of $X$. In general, this fails to be true. It is known that the fundamental group  of a locally complicated Peano continuum may contain nontrivial elements that are persistently undetectable by polyhedral approximations, at all scales. However, we show that such failure is not inherently local.
\end{abstract}

\maketitle

\vspace{-10pt}

\section{Introduction}

\noindent The process of approximating topological spaces with finite  polyhedra at progressively smaller scales is an important tool in many areas of mathematics,  ranging from shape theory to persistent homology. One goal is to calculate algebraic invariants across a sequence of approximating polyhedra in order to draw conclusions about the underlying space. Common obstacles include local geometric features with global impact. Therefore, such methods work best for path-connected, locally path-connected, compact metric spaces, i.e., Peano continua.

Still, the Griffiths twin cone, for example, has uncountable first homology, although it is a one-point union of two contractible Peano continua, and all of its fundamental group elements appear trivial in all polyhedral approximations. There is a well-known phenomenon at play here (see \S\ref{HPol}): If a fixed element $1\not=g\in\pi_1(X,\ast)$ of the fundamental group of a Peano continuum $X$ has the property that for every $\epsilon>0$, $g$ can be represented by a finite product of path-conjugates of $\epsilon$-small loops, then $g$ will not be picked up by any polyhedral approximation.

Indeed, there are many Peano continua in the literature whose fundamental groups contain nontrivial elements that are undetectable by polyhedral approximations. However, in all known examples, this is due to one of the following two local failures (cf.~\cite{BFi2020,FRVZ}):

\begin{itemize}
\item[(1)] {\em Point-local failure:} There is an $x\in X$ and a $1\not=g\in\pi_1(X,\ast)$ such that for every neighborhood $U$ of $x$, $g$ can be represented by a finite product of path-conjugates of loops that lie in $U$.
\item[(2)] {\em Path-local failure:} There is an essential loop $\alpha$ in $X$  such that for every $\epsilon>0$, there is an inessential loop $\alpha'$ $\epsilon$-close to $\alpha$.
\end{itemize}

We present an example of a two-dimensional Peano continuum $\mathbb{M}$ in $\mathbb{R}^3$, whose fundamental group contains uncountably many elements that are undetectable by polyhedral approximations (Theorem~\ref{noninj}), while satisfying neither (1) nor (2) (Corollary~\ref{main} and Theorem~\ref{T1}). Specifically, $\mathbb{M}$ is the mapping torus of the shift-map $f:\mathbb{D}\rightarrow \mathbb{D}$ on the dyadic arc space $\mathbb{D}$ (see \S\ref{Shift}).

A key step in our proof is Theorem~\ref{roundabout}, which is an adaption of a lemma used in the calculation of the first homology of the Griffiths twin cone in \cite{EF}. For a one-dimensional space $X$, it allows us to characterize the elements of the subgroups $\pi(U)\leqslant \pi_1(X,\ast)$, generated by path-conjugates of loops that lie in  $U\subseteq X$, in terms of their reduced representatives (Corollary~\ref{characterization}).

We close with a brief discussion of and comparison with other local properties of fundamental groups (see \S\ref{compare}).

\section{Polyhedral approximations\label{HPol}}

\noindent Let $X$ be a Peano continuum, that is, a path-connected, locally path-connected, compact metric space. Fix a basepoint $x_0\in X$. Then we can find (and fix) a sequence ${\mathcal U}_1, {\mathcal U}_2, {\mathcal U}_3, \dots$  of finite open covers of $X$ such that
\begin{itemize}
\item every $U\in {\mathcal U}_i$ is path connected (and not empty);
\item every $U\in {\mathcal U}_i$ has $diam(U)<\frac{1}{i}$;
\item ${\mathcal U}_{i+1}$ refines ${\mathcal U}_i$, i.e., $\forall\; V\in {\mathcal U}_{i+1}$ $\exists\; U\in {\mathcal U}_i:$ $V\subseteq U$;
\item for every $x\in X$, there is a designated $U^i_x\in {\mathcal U}_i$ with $x\in U^i_x$;
\item $U^i_{x_0}$ is the only element of ${\mathcal U}_i$ that contains $x_0$.
\end{itemize}

Let $P_i$ be the geometric realization of the finite abstract simplicial complex whose vertices are the elements of ${\mathcal U}_i$ and whose simplices are spanned by the (nonempty) subsets $\{U_0,U_1,\dots,U_n\}\subseteq {\mathcal U}_i$ with $U_0\cap U_1\cap\cdots \cap U_n\not=\emptyset$. The vertex $p_i=U^i_{x_0}$ will serve as the basepoint of $P_i$.

Choose functions $\psi_i:(P_{i+1},p_{i+1})\rightarrow (P_i,p_i)$ by first mapping each vertex $V\in {\mathcal U}_{i+1}$ to  some vertex  $U\in {\mathcal U}_i$ with $V\subseteq U$ and then extending the assignment linearly on all simplices of $P_{i+1}$. The induced homomorphism $$\psi_{i\#}:\pi_1(P_{i+1},p_{i+1})\rightarrow \pi_1(P_i,p_i)$$ between fundamental groups is independent of the choice of $\psi_{i}$ \cite{MS}.
The inverse limit \[\check{\pi}_1(X,x_0)=\lim_{\longleftarrow} \left( \pi_1(P_1, p_1) \stackrel{\psi_{1\#}}{\longleftarrow} \pi_1(P_2, p_2) \stackrel{\psi_{2\#}}{\longleftarrow} \pi_1(P_3, p_3) \stackrel{\psi_{3\#}}{\longleftarrow}\cdots\right)\]
is called the {\em first \v{C}ech homotopy group} of $X$ at $x_0$ and is independent of all preceding choices \cite{MS}.

Define $\varphi_i:X\rightarrow P_i$ based on any partition of unity $\{f_U\}_{U\in {\mathcal U}_i}$ subordinate to ${\mathcal U}_i$, by mapping $x\in X$ to the point of $P_i$ whose barycentric coordinate with respect to each vertex $U$ equals $f_U(x)$. In particular, $\varphi_i(x_0)=p_i$. Then the homomorphism $$\varphi_{i\#}:\pi_1(X,x_0)\rightarrow \pi_1(P_i, p_i)$$ is independent of the chosen partition of unity and $\psi_{i\#}\circ \varphi_{i+1\#}=\varphi_{i\#}$ \cite{MS}. Thus, we obtain a canonical homomorphism \[\Phi_{X,x_0}=(\varphi_{i\#})_i:\pi_1(X,x_0)\rightarrow \check{\pi}_1(X,x_0).\]

\begin{definition}
We call $\left\{ {\mathcal U}_i, U^i_x,P_i,p_i, \psi_i, \varphi_i \right\}$ a {\em nice polyhedral expansion} of $(X,x_0)$. (It is cofinal in  the more general \v{C}ech expansion \cite{MS}.)
\end{definition}

\begin{proposition}\cite{K} Suppose $X$ is a semilocally simply-connected Peano continuum and $x_0\in X$. Let $\left\{ {\mathcal U}_i, U^i_x,P_i,p_i, \psi_i, \varphi_i \right\}$ be a nice polyhedral expansion of $(X,x_0)$. Then there is an $n\in \mathbb{N}$ such that for every $i\geqslant n$, the homomorphism $\varphi_{i\#}:\pi_1(X,x_0)\rightarrow \pi_1(P_i, p_i)$ is an isomorphism. In particular, $\Phi_{X,x_0}:\pi_1(X,x_0)\rightarrow \check{\pi}_1(X,x_0)$ is an isomorphism.
\end{proposition}

\begin{proof}[Idea of Proof]
Choose $n$ large enough so that for every $i\geqslant n$ and every $U,V\in {\mathcal U}_i$ with $U\cap V\not=\emptyset$, the inclusion induced homomorphism $incl_{\#}:\pi_1(U\cup V,\ast)\rightarrow \pi_1(X,\ast)$ is trivial (cf. \cite[Theorem 4.1]{FZ2007}).
\end{proof}

\begin{remark} $\ker \Phi_{X,x_0}=1$ if $X$ is one-dimensional \cite{EK} or planar \cite{FZ2005}.
\end{remark}

\begin{remark} \label{Griffiths} There are many  examples of Peano continua $X$ for which $\ker \Phi_{X,x_0}\not=1$. (An extreme example is the Griffiths twin cone, where $\ker \Phi_{X,x_0}=\pi_1(X,x_0)$ \cite[2.5.18]{S} with $H_1(X)=\mathbb{Z}^\mathbb{N}/\bigoplus_\mathbb{N} \mathbb{Z}$ \cite{EF}.)
 Given a nice polyhedral expansion $\left\{ {\mathcal U}_i,  U^i_x,P_i,p_i, \psi_i, \varphi_i \right\}$ of $(X,x_0)$, note that for $g\in \pi_1(X,x_0)$, we have $g\in \ker \Phi_{X,x_0}$ if and only if $\varphi_{i\#}(g)=1\in \pi_1(P_i,p_i)$ for all $i\in \mathbb{N}$. Therefore, in these examples,  $\pi_1(X,x_0)$ contains nontrivial elements that cannot be detected by any of the polyhedral approximations. (To underscore the subtleties involved with the presence of such elements, we note that the fundamental group of the Griffiths twin cone is isomorphic to that of the triple cone, while no isomorphism between them is induced by a continuous map \cite{C}.)
\end{remark}

\begin{figure}[h!]
\includegraphics{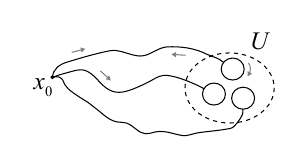}
\vspace{-23pt}

\caption{A representative of an  element of $\pi(U)$\label{pi-U}}
\end{figure}

\begin{definition}
Let $X$ be a path-connected space and $x_0\in X$.
For a path-connected subset $U\subseteq X$, we let $\pi(U)$ denote the (normal) subgroup of $\pi_1(X,x_0)$ which is generated by all elements $[\alpha\cdot \beta \cdot \alpha^{-}]$ represented by the concatenation of some path $\alpha:([0,1],0)\rightarrow (X,x_0)$ with $\alpha(1)\in U$, some loop $\beta:([0,1],\{0,1\})\rightarrow (U,\alpha(1))$, and the reverse path $\alpha^{-}(t)=\alpha(1-t)$. (See Figure~\ref{pi-U}.)
\end{definition}
\vspace{-25pt}

\vspace{-7pt}

\begin{proposition}\label{ker}\cite[Theorem 6.1]{BFa} Let $X$ be a Peano continuum and $x_0\in X$.
Let $\left\{ {\mathcal U}_i, U^i_x,P_i,p_i, \psi_i, \varphi_i \right\}$ be a nice polyhedral expansion of $(X,x_0)$. Say, ${\mathcal U}_i=\{U^i_0, U^i_1, U^i_2,\dots, U^i_{n_i}\}$.
Then \[\ker \Phi_{X,x_0} =\bigcap_{i\in \mathbb{N}} \pi(U^i_0)\pi( U^i_1)\pi(U^i_2)\cdots \pi(U^i_{n_i}). \]
\end{proposition}

For related constructions used in persistent homology, see \cite{Virk}. The proof of the following lemma is an elementary exercise.

\begin{lemma}\label{loops} Let $g\in \pi_1(X,x_0)$ and let $U_0, U_1, \dots, U_n\subseteq X$ be path-connected subsets.
Then
$g\in \pi(U_0)\pi( U_1)\pi(U_2)\cdots \pi(U_n)$ if and only if there are paths $\alpha_i$ and loops $\beta_i$ with $g=[\alpha_1\cdot \beta_1\cdot \alpha_2\cdot \beta_2\cdots \alpha_k\cdot \beta_k \cdot \alpha_{k+1}]$ such that each $\beta_i$ lies in some $U_j$ and $[\alpha_1\cdot \alpha_2\cdots \alpha_{k+1}]=1$.
\end{lemma}

\begin{remark} \label{eps-T1} Let $(X,d)$ be a path-connected and locally path-connected metric space.
Let $[\alpha]\in \pi_1(X,x_0)$ be such that for every $\epsilon>0$, there is a null-homotopic loop  $\alpha':([0,1],\{0,1\})\rightarrow (X,x_0)$ with $d(\alpha(t),\alpha'(t))<\epsilon$ for all $t\in [0,1]$.
Let $U_0, U_1, \dots, U_n$ be  path-connected open subsets of $X$ with $\alpha([0,1])\subseteq \bigcup_{j=0}^n U_j$.
Then there is a null-homotopic loop $\alpha'$, a subdivision $0=t_0<t_1<\dots <t_k=1$ of $[0,1]$, and loops $\beta_i$ with $[\alpha]=[\beta_1\cdot \alpha'_1\cdot \beta_2\cdot \alpha'_2\cdots \beta_k\cdot \alpha'_k]$ where each $\beta_i$ lies in some $U_j$ and $\alpha'_i=\alpha'|_{[t_{i-1},t_i]}$ for $i=1,2,\dots, k$. In particular, $[\alpha]\in \pi(U_0)\pi( U_1)\pi(U_2)\cdots \pi(U_n)$.

As in \cite{BFa2015}, one shows that the following two statements are equivalent (even if $X$ is not locally path connected):
\begin{itemize}
\item[(a)] There exists $1\not=[\alpha]\in \pi_1(X,x_0)$, such that for every $\epsilon>0$, there is a null-homotopic loop  $\alpha':([0,1],\{0,1\})\rightarrow (X,x_0)$ with $d(\alpha(t),\alpha'(t))<\epsilon$ for all $t\in [0,1]$.
\item[(b)] $\pi_1(X,x_0)$ fails to be T$_1$ in the quotient topology of the compact-open topology on the loop space $\Omega(X,x_0)$.
\end{itemize}
In particular, if $\ker \Phi_{X,x_0}=1$, then $\pi_1(X,x_0)$ is T$_1$. (See also \cite{FRVZ}.)
\end{remark}

\begin{remark}\label{local}
For every nice polyhedral expansion $\left\{ {\mathcal U}_i,  U^i_x,P_i,p_i, \psi_i, \varphi_i \right\}$ of a Peano continuum $(X,x_0)$ and for every $x\in X$, we clearly have \[\bigcap_{i\in \mathbb{N}} \pi(U^i_x)\leqslant \bigcap_{i\in \mathbb{N}} \pi(U^i_0)\pi( U^i_1)\pi(U^i_2)\cdots \pi(U^i_{n_i}),\] where ${\mathcal U}_i=\{U^i_0, U^i_1, U^i_2,\dots, U^i_{n_i}\}$.
\end{remark}

Considering Proposition~\ref{ker} and Lemma~\ref{loops}, one may be led to believe that the property ``$\ker \Phi_{X,x_0}\not=1$'' is a local property of  $X$. However, in view of Remark~\ref{eps-T1} and Remark~\ref{local}, one may also wonder: How local is this property, really?
Specifically, our goal is to answer the following question.

\begin{question}\label{Q} Is there a Peano continuum $X$ such that
\begin{itemize}\setlength{\itemsep}{5pt}
\item $\ker \Phi_{X,x_0}\not=1$, and
\item $\bigcap_{i\in \mathbb{N}} \pi(U^i_x)=1$  for all $x\in X$, and
\item $\pi_1(X,x_0)$ is T$_1$?
\end{itemize}
\end{question}

We give a positive answer in Theorem~\ref{noninj}, Corollary~\ref{main}, and Theorem~\ref{T1}.

\section{\label{Shift}The shift-map $f:\mathbb{D}\rightarrow \mathbb{D}$ on the dyadic arc space $\mathbb{D}$}

\begin{figure}[h!]
\includegraphics[scale=3]{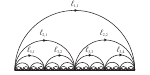}
\caption{The space $\mathbb{D}$ and paths $\ell_{n,j}$ \label{dyadic}}
\end{figure}

\noindent For each $n\in \mathbb{N}$ and $j\in \{1,2,\dots, 2^{n-1}\}$, we define \[\mathbb{D}_{n,j}=\{(x,y)\in \mathbb{R}^2\mid \left(x-\frac{2j-1}{2^n}\right)^2+y^2=\left(\frac{1}{2^n}\right)^2,\;  y\geqslant 0\}.\]
Let $B=[0,1]\times \{0\}\subseteq \mathbb{R}^2$ and $\mathbb{D}=B\cup \bigcup_{n,j} \mathbb{D}_{n,j}\subseteq \mathbb{R}^2$ with basepoint $d_0=(0,0)\in \mathbb{D}$.
Let $\ell_{n,j}:[0,1]\rightarrow \mathbb{D}_{n,j}$ be the homeomorphism given by $\ell_{n,j}(t)=(\frac{j-1+t}{2^{n-1}},\frac{\sqrt{t-t^2}}{2^{n-1}})$. (See Figure~\ref{dyadic}.)

We define the {\em shift map} $f:\mathbb{D}\rightarrow \mathbb{D}$ by $f(x,0)=(x,0)$ for $(x,0)\in B$ and $f\circ \ell_{n,j}=\ell_{n+1,2j-1}\cdot \ell_{n+1,2j}$ for all $n\in \mathbb{N}$ and $j\in\{1,2,\dots,2^{n-1}\}$. That is, for $(x,y)\in \mathbb{D}_{n,j}$, we have $f(x,y)=(x,y')\in \mathbb{D}_{n+1,2j-1}\cup \mathbb{D}_{n+1,2j}$.

\begin{lemma}\label{f-inj}
$f_\#:\pi_1(\mathbb{D},d_0)\rightarrow \pi_1(\mathbb{D},d_0)$ is injective.
\end{lemma}

\begin{proof}
Let $E_n$ denote the finite graph consisting of the union of $B$ and all $\mathbb{D}_{k,j}$ with $1\leqslant k \leqslant n$ and $j\in\{1,2,\dots,2^{k-1}\}$. Define retractions $r_n:\mathbb{D}\rightarrow E_n$ by $r_n(x,y)=(x,0)\in B$ for $(x,y)\in \mathbb{D}_{k,j}$ with $k>n$ and $j\in\{1,2,\dots,2^{k-1}\}$.
Then $\mathbb{D}$ is the inverse limit of $ E_1 \stackrel{r_1|_{E_2}}{\longleftarrow} E_2\stackrel{r_2|_{E_3}}{\longleftarrow}  \cdots$ with projections $r_n:\mathbb{D}\rightarrow E_n$. It is a straightforward exercise to show that each $(f|_{E_n})_\#:\pi_1(E_n,d_0)\rightarrow \pi_1(E_{n+1},d_0)$ is an injective homomorphism of free groups, inducing an injective homomorphism $\displaystyle \lim_{\longleftarrow} ((f|_{E_n})_\#)_n$ between the inverse limits in the following commutative diagram:
\[
\xymatrix{
\pi_1(\mathbb{D},d_0) \ar[d]_{f_\#} \ar[rr]^(.25){(r_{n\#})_n} && \displaystyle \lim_{\longleftarrow} \left(\pi_1(E_1,d_0) \stackrel{(r_1|_{E_2})_\#}{\longleftarrow} \pi_1(E_2,d_0)\stackrel{(r_2|_{E_3})_\#}{\longleftarrow}  \cdots\right) \ar[d]^{\displaystyle \lim_{\longleftarrow} ((f|_{E_n})_\#)_n} \\
\pi_1(\mathbb{D},d_0) \ar[rr]^(.25){(r_{n+1\#})_n} && \displaystyle \lim_{\longleftarrow} \left(\pi_1(E_2,d_0) \stackrel{(r_2|_{E_3})_\#}{\longleftarrow} \pi_1(E_3,d_0)\stackrel{(r_3|_{E_4})_\#}{\longleftarrow}  \cdots\right)
}
\]
Since $\mathbb{D}$ is  one-dimensional, $(r_{n\#})_n$ is also injective~\cite{EK}. Consequently, $f_\#$ is injective.
\end{proof}

Since all fibers of $f$ are totally path disconnected, we observe:

\begin{lemma} A path $\alpha$ in $\mathbb{D}$ is constant if and only if $f\circ \alpha$ is constant.
\end{lemma}

\begin{lemma}\label{uniquelift}
Let $\alpha:[0,1]\rightarrow \mathbb{D}$ be a nonconstant path. If there are two paths $\beta_1, \beta_2:[0,1]\rightarrow \mathbb{D}$ with $f\circ \beta_i=\alpha$ for $i=1,2$, then $\beta_1=\beta_2$.
\end{lemma}

\begin{proof}
Suppose there is a $t\in [0,1]$ such that $\beta_1(t)\not=\beta_2(t)$. Since $f\circ \beta_1(t)=f\circ \beta_2(t)$, we may assume that $\beta_1(t)=(\frac{2j-1}{2^{n}},\frac{1}{2^{n}})$ and $\beta_2(t)=(\frac{2j-1}{2^{n}},0)$ for some $n\in \mathbb{N}$ and $j\in \{1,2,\dots,2^{n-1}\}$. Since $\beta_1$ is not constant, there is an $s\in [0,1]$ such that $\beta_1(s)\not=\beta_1(t)$. We may assume that $s>t$. Choose $m\in \mathbb{N}$ such that  $||\beta_i(s_1)-\beta_i(s_2)||<\frac{1}{2^{n+1}}$ for all $s_1, s_2\in [0,1]$ with $|s_1-s_2|\leqslant 1/m$ and $i\in \{1,2\}$. Choose $s_1, s_2\in [t,s]$ with $|s_1-s_2|\leqslant 1/m$ such that $\beta_1|_{[t,s_1]}$ is constant and $\beta_1(s_2)\not=\beta_1(t)$. Then $\beta_2|_{[t,s_1]}$ is also constant. Therefore, $f\circ \beta_1(s_2)\in \left(\mathbb{D}_{n+1,2j-1}\cup \mathbb{D}_{n+1,2j}\right)\setminus\{(\frac{2j-1}{2^n},0)\}$, while $f\circ\beta_2(s_2)$ lies in closed region of $\mathbb{R}^2$ bounded by the three arcs $\mathbb{D}_{n+2,4j-2}$, $\mathbb{D}_{n+2,4j-1}$, and $[\frac{4j-3}{2^{n+1}},\frac{4j-1}{2^{n+1}}]\times\{0\}$. Hence, $f\circ\beta_1(s_2)\not=f\circ\beta_2(s_2)$; a contradiction.
\end{proof}

\begin{definition} A path $\alpha:[a,b]\rightarrow X$ is called {\em reduced} if $\alpha$ is either constant or if for all $a\leqslant s<t\leqslant b$ with $\alpha(s)=\alpha(t)$, the loop $\alpha|_{[s,t]}$ is not null-homotopic. For two paths $\alpha:[a,b]\rightarrow X$ and $\beta:[c,d]\rightarrow X$, we write $\alpha\equiv \beta$ (and call the paths {\em equivalent})
 if there is an increasing homeomorphism $\gamma:[a,b]\rightarrow [c,d]$ such that $\alpha=\beta\circ \gamma$.
\end{definition}

\begin{remark}\label{letter-in-letter} If $\alpha:[a_0,b_0]\rightarrow X$ is a non-constant reduced path in a metric space $X$ and  $[a_1,b_1]\subseteq [a_0,b_0]$ with $\alpha\equiv \alpha|_{[a_1,b_1]}$, then $[a_1,b_1]=[a_0,b_0]$. (Indeed, let $\gamma:[a_0,b_0]\rightarrow [a_1,b_1]$ be an increasing homeomorphism with $\alpha=\alpha|_{[a_1,b_1]}\circ \gamma$. Recursively, put $a_{i+1}=\gamma(a_i)$ and $b_{i+1}=\gamma(b_i)$ for $i\geqslant 1$. Then $a_0\leqslant a_1\leqslant a_2\leqslant \cdots \leqslant b_2\leqslant b_1\leqslant b_0$ with $\alpha|_{[a_{i},a_{i+1}]}\equiv \alpha|_{[a_{i-1},a_i]}$ and $\alpha|_{[b_{i+1},b_{i}]}\equiv \alpha|_{[b_{i},b_{i-1}]}$. Since $(a_i)_i$ converges to some $a$, the continuity of $\alpha$ implies that $\alpha|_{[a_0,a_1]}$ is constant. Hence, as $\alpha$ is reduced, $a_1=a_0$. Similarly, $b_1=b_0$.)
\end{remark}

\begin{remark} We recall from \cite{E} that in a one-dimensional metric space, every path is path-homotopic (within its own image) to a reduced path, which is unique up to equivalence.
\end{remark}

\begin{definition} A metric space $X$ is called a {\em geodesic space} if every two distinct points of $X$ are the endpoints of some geodesic (a subset of $X$ that is isometric to some closed interval $[a,b]\subseteq \mathbb{R}$).

A topological space $X$ is called {\em uniquely arcwise connected} if every two distinct points of $X$ are the endpoints of a unique arc (a subset of $X$ that is homeomorphic to the unit interval $[0,1]\subseteq \mathbb{R}$).

A uniquely arcwise connected geodesic space is called an {\em $\mathbb{R}$-tree}.
\end{definition}

\begin{remark} $\mathbb{R}$-trees are locally path connected and contractible.
\end{remark}
\begin{remark}\label{tree}
We recall from \cite{FZ2007} that for every path-connected one-dimensional separable metric space $X$, there is an $\mathbb{R}$-tree $\widetilde{X}$ and a continuous surjection $p:\widetilde{X}\rightarrow X$ with the following property:

For every path-connected and locally path-connected
space $Y$, for every continuous function $g : (Y, y) \rightarrow (X, x)$ with
$g_\#(\pi_1(Y, y)) = 1$, and for every $\widetilde{x}\in \widetilde{X}$ with $p(\widetilde{x}) = x$, there is
a unique continuous lift $\widetilde{g} : (Y, y) \rightarrow  (\widetilde{X} , \widetilde{x})$ with $p \circ \widetilde{g} = g$.

In particular, null-homotopic loops lift to (null-homotopic) loops. Moreover, the automorphism group Aut$(\widetilde{X}\stackrel{p}{\rightarrow} X)$, i.e., the group of homeomorhisms $h:\widetilde{X}\rightarrow \widetilde{X}$ with $p\circ h=p$, is isomorphic to $\pi_1(X,\ast)$.
\end{remark}

\begin{lemma}\label{reduced} Let $\alpha$ be a loop in $\mathbb{D}$ based at $d_0$.
Then  $\alpha$ is reduced if and only if $f\circ \alpha$ is reduced.
\end{lemma}

\begin{proof}
If $\alpha$ is not reduced, then neither is $f\circ \alpha$. Now assume that $\alpha$ is non-constant and reduced. Let $A=\alpha^{-1}(B)$. Then for every component $(a,b)$ of $[0,1]\setminus A$, either $\alpha|_{[a,b]}\equiv \ell_{n,j}$ or $\alpha|_{[a,b]}\equiv \ell_{n,j}^{-}$ for some $n$ and $j$,
 so that either $f\circ \alpha|_{[a,b]}\equiv \ell_{n+1,2j-1}\cdot \ell_{n+1,2j}$ or $f\circ \alpha|_{[a,b]}\equiv  \ell_{n+1,2j}^-\cdot \ell_{n+1,2j-1}^-$.
Suppose, to the contrary, that there are $0\leqslant s<t\leqslant 1$ such that $f\circ\alpha|_{[s,t]}$ is a null-homotopic loop in $\mathbb{D}$. If $\alpha(s)=\alpha(t)$, then $\alpha|_{[s,t]}$ is a loop in $\mathbb{D}$ that is not nullhomotopic, contradicting Lemma~\ref{f-inj}.  So, $\alpha(s)\not=\alpha(t)$. Since $f\circ\alpha(s)=f\circ\alpha(t)$, we may assume that $\alpha(s)=(\frac{2j-1}{2^{n}},\frac{1}{2^{n}})$ and $\alpha(t)=(\frac{2j-1}{2^{n}},0)$ for some $n$ and $j$. Let $(a,b)$ be the component of $[0,1]\setminus A$ with $a<s<b<t$. We may assume that $\alpha|_{[a,b]}\equiv \ell_{n,j}$, so that $f\circ\alpha|_{[s,b]}\equiv \ell_{n+1,2j}$.

Let $p:\widetilde{\mathbb{D}}\rightarrow \mathbb{D}$ be as in Remark~\ref{tree}.
Consider $g=f\circ \alpha|_{[s,t]}:[s,t]\rightarrow \mathbb{D}$ and choose $\widetilde{x}\in \widetilde{\mathbb{D}}$ with $p(\widetilde{x})=g(s)$. Let $\widetilde{g}:([s,t],s)\rightarrow (\widetilde{\mathbb{D}},\widetilde{x})$ be the lift with $p\circ\widetilde{g}=g$. Since $g$ is a null-homotopic loop in $\mathbb{D}$, we have $\widetilde{g}(s)=\widetilde{g}(t)$. Then $\widetilde{\ell}_{n+1,2j}=\widetilde{g}|_{[s,b]}$ is an embedding onto an arc in $\widetilde{\mathbb{D}}$. Since $\widetilde{\mathbb{D}}$ is uniquely arcwise connected, there are $b\leqslant a'<t'\leqslant t$ such that $\widetilde{g}(a')=\widetilde{g}(b)$, $\widetilde{g}(t')=\widetilde{g}(s)$, and $\widetilde{g}((a',t'))\subseteq \widetilde{g}((s,b))$.
Note that $p\circ \widetilde{g}|_{[a',t']}=f\circ \alpha|_{[a',t']}:[a',t']\rightarrow p(\widetilde{g}([s,b]))=\mathbb{D}_{n+1,2j}$.
By definition of $f$ (and because $\alpha$ is reduced) there is a $b'\in(t',1]$ such that $\alpha|_{[a',b']}\equiv \ell_{n,j}^-$ and
$\widetilde{g}|_{[a',t']}\equiv\widetilde{\ell}_{n+1,2j}^-$. In particular, $\alpha(t')=\alpha(s)$.
However,
since $\widetilde{\mathbb{D}}$ is simply connected, $\widetilde{g}|_{[s,t']}$ is a null-homotopic loop in $\widetilde{\mathbb{D}}$. Hence, $p\circ\widetilde{g}|_{[s,t']}=f\circ \alpha|_{[s,t']}$ is a null-homotopic loop in $\mathbb{D}$. As above, this contradicts Lemma~\ref{f-inj}.
\end{proof}

\section{A characterization of $\pi(U)$ for one-dimensional spaces}\label{key}

\begin{definition}\label{def-cancel} Let $X$ be a one-dimensional metric space with basepoint $x_0\in X$. Let $U\subseteq X$. Let $\alpha:([a,b],\{a,b\})\rightarrow (X,x_0)$ be a loop.
A {\em cancellation of $\alpha$ relative to $U$} consists of
\begin{itemize}
\item a subdivision $a=t_0<t_1<\cdots<t_n=b$,
\item a partition $\{1,2,\dots,n\}=C\cup \overline{C}\cup R$,
\item a bijective function $\varphi:C\rightarrow \overline{C}$, and
\item for each $i\in C$, an arc $S_i$ (which we may choose to be a semicircle) in the upper half-plane $\mathbb{R}\times[0,\infty)$
with
$\partial S_i= S_i\cap \left(\mathbb{R}\times\{0\}\right)$, $\partial S_i\cap \left((t_{i-1},t_i)\times\{0\}\right)\not= \emptyset$ and  $\partial S_i\cap \left((t_{\varphi(i)-1},t_{\varphi(i)})\times\{0\}\right)\not=\emptyset$,

\end{itemize}

such that

\begin{itemize}
\item[(1)] for all $1\leqslant i\leqslant n$,  $\alpha_i=\alpha|_{[t_{i-1},t_i]}$ is reduced,
\item[(2)] for all $i\in R$, Im$(\alpha_i)\subseteq U$,
\item[(3)] for all $i\in C$,  $\alpha_i\equiv \alpha_{\varphi(i)}^{-}$, and
\item[(4)] for all $i,j\in C$ with $i\not=j$, $S_i\cap S_j=\emptyset$.
\end{itemize}
\end{definition}

\pagebreak

\begin{figure}[h!]
 \includegraphics[scale=0.9]{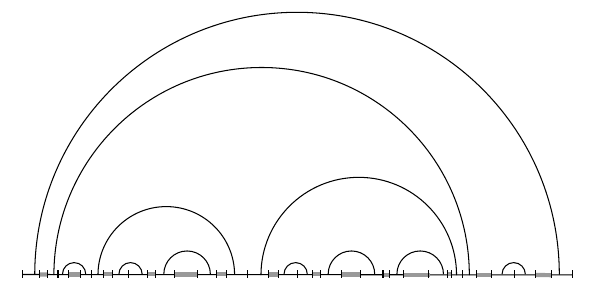} \hspace{-30pt} \includegraphics[scale=0.8]{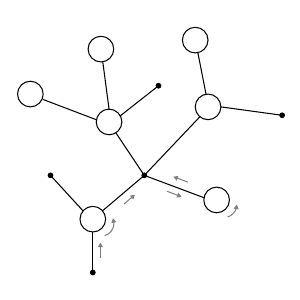}

\caption{\label{roundabout-fig} Intervals that are paired by semicircles are identified in reverse orientation. All other (gray) intervals form the roundabouts that map into $U$.}
\end{figure}

\begin{remark}\label{representation}
A loop that has a cancellation relative to $U$ factors through a loop going counterclockwise around a planar simplicial tree, some of whose vertices have been ``blown up'' to circles (like roundabouts) that map into $U$. The semicircles $S_i$ prescribe a cancellation pattern for the edge-path in the underlying tree. (See Figure~\ref{roundabout-fig}.)
\end{remark}

\begin{theorem}\label{roundabout}
Let $X$ be a path-connected one-dimensional separable metric space with basepoint $x_0\in X$ and $U\subseteq X$. Suppose $1\not=[\alpha]\in \pi_1(X,x_0)$. If there is a cancellation of $\alpha$ relative to $U$, then there is a cancellation of the reduced representative for $[\alpha]$ relative to $U$.
\end{theorem}

\begin{proof} Let $\alpha:([a,b],\{a,b\})\rightarrow (X,x_0)$ be an essential loop.
Suppose there is a cancellation of $\alpha$ relative to $U$, as described in Definition~\ref{def-cancel}, with $a=t_0<t_1<\cdots<t_n=b$, $\alpha_i=\alpha|_{[t_{i-1},t_i]}$, $\{1,2,\dots,n\}=C\cup \overline{C}\cup R$, $\varphi:C\rightarrow \overline{C}$, and $S_i$ ($i\in C$). We may assume that each $\alpha_i$ is non-constant. Let $p:\widetilde{X}\rightarrow X$ be as in Remark~\ref{tree}. Fix any lift $\widetilde{\alpha}:[a,b]\rightarrow \widetilde{X}$ with $p\circ \widetilde{\alpha}=\alpha$. Put $\widetilde{\alpha}_i=\widetilde{\alpha}|_{[t_{i-1},t_i]}$ so that $\alpha_i=p\circ \widetilde{\alpha}_i$.
 Since each $\alpha_i$ is reduced and since $\widetilde{X}$ is simply connected, it follows that each $\widetilde{\alpha}_i$ is an embedding onto an arc. Reparametrizing $\alpha$ as necessary, we may assume that each $\widetilde{\alpha}_i$ is (an isometric embedding onto) a geodesic. Let $\widetilde{\beta}:[0,d]\rightarrow \widetilde{X}$ be the unique geodesic with $\widetilde{\beta}(0)=\widetilde{\alpha}(a)$ and $\widetilde{\beta}(d)=\widetilde{\alpha}(b)$. Then $\beta=p\circ \widetilde{\beta}$ is a reduced representative of $[\alpha]$.

 We wish to show that there is a cancellation of $\beta$ relative to $U$.
The argument is the same as in the proof of \cite[Lemma~4.3 ]{EF}. However, one needs to check that the transformations described there, do not only respect Properties (1)--(3), but also Property~(4). We sketch the idea:

Let $r$ denote the number of indices $i\in \{1,2,\dots, n-1\}$ for which $\widetilde{\alpha}|_{[t_{i-1},t_{i+1}]}$ is not a geodesic. If $r=0$, then $\widetilde{\alpha}$ is a geodesic and $\beta\equiv\alpha$, in which case we are done. Otherwise, we recursively subject $\widetilde{\alpha}$ to the following transformation, which replaces $\widetilde{\alpha}$ with $\widetilde{\alpha}|_{[a,a_1]\cup [b_1,b]}$ (where we abuse notation) for some $\widetilde{\alpha}|_{[a_1, t_i]}\equiv \widetilde{\alpha}|_{[t_i,b_1]}^-$, while retaining a cancellation relative to $U$ for the projection into $X$ (which is homotopic to $\alpha$) on a refined subdivision (which includes $a_1$ and $b_1$), eventually reducing the pair $(r,n)$ in the lexicographical ordering.
\vspace{10pt}

{\bf Step~1.}
Let $i$ be the largest index for which $\widetilde{\alpha}|_{[t_{i-1},t_{i+1}]}$ is not a geodesic. Since $\widetilde{X}$ is uniquely arcwise connected, there are unique points  $a_1\in [t_{i-1},t_i)$ and $b_1\in (t_i,b]$ such that the arc $[\widetilde{\alpha}(t_{i-1}),\widetilde{\alpha}(b)]$ from $\widetilde{\alpha}(t_{i-1})$ to $\widetilde{\alpha}(b)$ in $\widetilde{X}$ equals $[\widetilde{\alpha}(t_{i-1}),\widetilde{\alpha}(a_1)]\cup[\widetilde{\alpha}(b_1), \widetilde{\alpha}(b)]$. In particular, $\widetilde{\alpha}(a_1)=\widetilde{\alpha}(b_1)$ and $\widetilde{\alpha}|_{[a_1,t_i]}\equiv \widetilde{\alpha}|_{[t_i,b_1]}^-$. Say, $b_1\in (t_{j-1},t_j]$. Define $\alpha_{(i,1)}=\alpha|_{[t_{i-1},a_1]}$, $\alpha_{(i,2)}=\alpha_{[a_1,t_i]}$, $\alpha_{(j,1)}=\alpha|_{[t_{j-1},b_1]}$, and $\alpha_{(j,2)}=\alpha|_{[b_1,t_j]}$. We will use the notational shorthand $\alpha_i=\alpha_{(i,1)}\alpha_{(i,2)}$ and $\alpha_j=\alpha_{(j,1)}\alpha_{(j,2)}$. Then $\alpha_{(i,2)}\equiv(\alpha_{i+1}\alpha_{i+2}\cdots \alpha_{j-1}\alpha_{(j,1)})^-$. (See Figure~\ref{Step1}.) If $a_1=t_{i-1}$, then $\alpha_{(i,1)}$ is degenerate and we treat it as empty. Likewise, if $b_1=t_j$, we treat $\alpha_{(j,2)}$ as empty.

\begin{figure}[h!]
\hspace{-10pt} \includegraphics{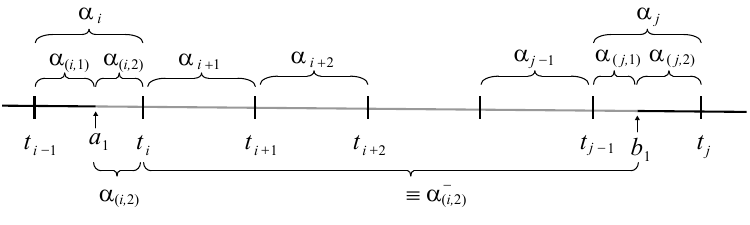}
\vspace{-15pt}

\caption{\label{Step1} The points $a_1$ and $b_1$ of Step 1}
\end{figure}

{\bf Step~2.} Suppose $i\in C$. Then $\alpha_{\varphi(i)}\equiv \alpha_i^-=\alpha_{(i,2)}^-\alpha_{(i,1)}^-$, so that there is a subdivision $\eta=\{t_{\varphi(i)-1}<t^\ast_{i+1}<t^\ast_{i+2}<\cdots<t^\ast_{j-1}<t^\ast_j\leqslant t_{\varphi(i)}\}$ such that $\alpha_{\varphi(i)}|_{[t_{\varphi(i)-1},t^\ast_{i+1}]}\equiv \alpha_{i+1}$, $\alpha_{\varphi(i)}|_{[t^\ast_{k-1},t^\ast_k]}\equiv \alpha_k$ for $i+2\leqslant k\leqslant j-1$, $\alpha_{\varphi(i)}|_{[t^\ast_{j-1},t^\ast_j]}\equiv \alpha_{(j,1)}$, and $\alpha_{\varphi(i)}|_{[t^\ast_j,t_{\varphi(i)}]}\equiv \alpha_{(i,1)}^-$.
Note that by Remark~\ref{letter-in-letter}, either $\varphi(i)<i$ or $\varphi(i)\geqslant j$.
Suppose $j\in \overline{C}$. Then $\alpha_{\varphi^{-1}(j)}\equiv \alpha_j^-=\alpha_{(j,2)}^-\alpha_{(j,1)}^-$, so that there is a $t_j^{\ast\ast}\in [t_{\varphi^{-1}(j)-1},t_{\varphi^{-1}(j)})$ such that $\alpha_{\varphi^{-1}(j)}|_{[t_{\varphi^{-1}(j)-1},t^{\ast\ast}_j]}\equiv \alpha_{(j,2)}^-$ and
$\alpha_{\varphi^{-1}(j)}|_{[t^{\ast\ast}_j,t_{\varphi^{-1}(j)}]}\equiv \alpha_{(j,1)}^-$.
Analogous statements hold if $i\in\overline{C}$ or $j\in C$.

Note that if $i\in R$, then $\alpha_i, \alpha_{i+1}, \dots, \alpha_{j-1}$ and $\alpha_{(j,1)}$ all lie in $U$, so that there is no need to
make replacement (iv) in Step~3.1 below.
 Likewise, if $j\in R$, then there is no need to make replacement (v) in Step~3.1 below.
 So, we may assume that $i\in C$ and $j\in \overline{C}$.
\vspace{10pt}

{\bf Step~3.1.} Suppose $\varphi(i)\not = j$.
\vspace{10pt}

We may reduce the domain of $\widetilde{\alpha}$ to $[a,a_1]\cup [b_1,b]$ and adjust the concatenation $\alpha_1\alpha_2\cdots\alpha_n$ to the new subdivision by carrying out the following steps, in this order (see Figure~\ref{Step3.1}):
\begin{itemize}
\item[(i)] replace $\alpha_i$ with $\alpha_{(i,1)}$ (or eliminate it, if $\alpha_{(i,1)}$ is empty)
\item[(ii)]  eliminate $\alpha_{i+1}, \alpha_{i+2}, \dots,\alpha_{j-1}$
\item[(iii)] replace $\alpha_j$ with $\alpha_{(j,2)}$ (or eliminate it, if $\alpha_{(j,2)}$ is empty)
\item[(iv)] replace $\alpha_{\varphi(i)}$ with $\widehat{\alpha}_{i+1}\widehat{\alpha}_{i+2}\cdots \widehat{\alpha}_{j-1}\widehat{\alpha}_{(j,1)}\widehat{\alpha}^-_{(i,1)}$, where $\widehat{\alpha}_{\bullet}\equiv \alpha_{\bullet}$
\item[(v)] If $\varphi^{-1}(j)\not\in \{i+1, i+2, \dots, j-1\}$, then replace $\alpha_{\varphi^{-1}(j)}$ with $\widehat{\alpha}^-_{(j,2)}\widehat{\alpha}^-_{(j,1)}$; otherwise, replace $\widehat{\alpha}_{\varphi^{-1}(j)}$ from step (iv) with
    $\widehat{\alpha}^-_{(j,2)}\widehat{\alpha}^-_{(j,1)}$
\end{itemize}
(Note that there is no $\widehat{\alpha}_{(i,2)}$.)

\begin{figure}[h!]
\parbox{5.5in}{
\hspace{-.5in} \includegraphics{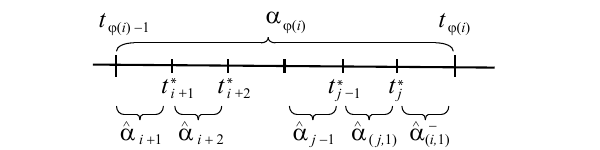} \hspace{-.45in} \includegraphics{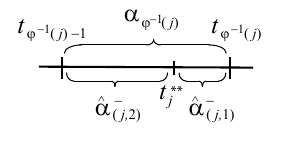}}
\caption{\label{Step3.1} The replacements (iv) and (v) of Step 3.1}
\end{figure}

 We now need to place and reposition some semicircles, without introducing intersections: (See Figure~\ref{slide}, where $\varphi(i)<i<j<\varphi^{-1}(j)$.)

\begin{itemize}
\item[(a)] We may assume that $S_i$ and $S_{\varphi^{-1}(j)}$ were already positioned so as to serve as $\widehat{S}_{(i,1)}$ and $\widehat{S}_{(j,2)}$, matching $\alpha_{(i,1)}$ with $\widehat{\alpha}^-_{(i,1)}$ and matching $\alpha_{(j,2)}$ with $\widehat{\alpha}^-_{(j,2)}$, respectively, if these semicircles are still needed.
\item[(b)] Slide the end of any semicircle that touches the domain of one of $\alpha_{i+1}, \alpha_{i+2}, \dots, \alpha_{j-1}$ along the $x$-axis and along
$S_i=\widehat{S}_{(i,1)}$ to the corresponding point of $\widehat{\alpha}_{i+1}, \widehat{\alpha}_{i+2}, \dots, \widehat{\alpha}_{j-1}$.

\item[(c)] Place a new semicircle $\widehat{S}_{(j,1)}$ (matching  $\widehat{\alpha}_{(j,1)}$ with $\widehat{\alpha}^-_{(j,1)}$) closely following the $x$-axis and the two semicircles $S_i=\widehat{S}_{(i,1)}$ and $S_{\varphi^{-1}(j)}=\widehat{S}_{(j,2)}$.
\end{itemize}
\vspace{-25pt}

\begin{figure}[h]
\includegraphics{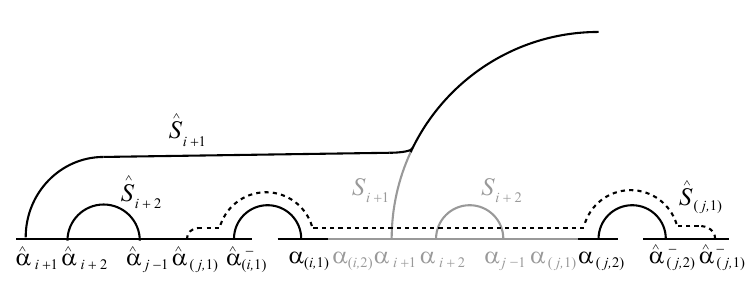}
\vspace{-30pt}

\caption{\label{slide} Reposition, slide and place semicircles}
\end{figure}

\pagebreak

{\bf Step~3.2.} Suppose $\varphi(i)= j$.
\vspace{10pt}

 For every $s\in C\cap \{i+1, i+2, \dots,j-1\}$, we have  $i+1\leqslant \varphi(s)\leqslant j-1$, because $S_s\cap S_i=\emptyset$. Likewise, if  $s\in \overline{C}\cap \{i+1, i+2, \dots,j-1\}$, then $i+1\leqslant \varphi^{-1}(s)\leqslant j-1$.

We now consider the common subdivision $\eta\cup\{b_1\}$ of $[t_{j-1},t_j]$. Since $[t_{\varphi(i)-1},t_{\varphi(i)}]=[t_{j-1},t_j]$ and $\alpha|_{[t^\ast_{j-1},t^\ast_j]}\equiv \alpha|_{[t_{j-1},b_1]}$, it follows from Remark~\ref{letter-in-letter} that $[t_{j-1},b_1]$ cannot properly contain $[t^\ast_{j-1},t^\ast_j]$.
Therefore, $b_1\leqslant t^\ast_j$. Moreover, if $b_1 = t^\ast_j$, then $j = i+1$ (so that there are no $\alpha_{i+1},\dots , \alpha_{j-1}$)
and $\alpha_i\equiv \alpha_{i+1}^-$, in which case we are done. Hence, we may assume that $b_1 < t^\ast_j$.
\vspace{10pt}

{\bf Case~1.} Suppose $b_1\leqslant t^\ast_{j-1}$. Let $k$ be minimal with $b_1\leqslant t^\ast_{k}$.
\vspace{20pt}

We may transfer the subdivision points $\eta\cap [t_{j-1},b_1]$ to $[b_1,t_j]$ along the correspondence of $\alpha_{(j,1)}$ with $\widehat{\alpha}_{(j,1)}$, where $\widehat{\alpha}_{\bullet}\equiv \alpha_{\bullet}$ as above. (See Figure~\ref{Step3.2Case1}.)
\begin{figure}[h!]
\hspace{.5in} \includegraphics{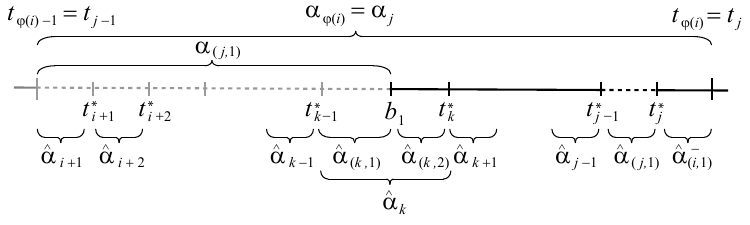}

\caption{\label{Step3.2Case1}Transfer of subdivision points (Step~3.2, Case1)}
\end{figure}

We  reduce the domain of $\widetilde{\alpha}$ to $[a,a_1]\cup [b_1,b]$ and adjust the concatenation $\alpha_1\alpha_2\cdots\alpha_n$ to the new subdivision by carrying out the following steps, in this order:
\begin{itemize}
\item[(i)] replace $\alpha_i$ with $\alpha_{(i,1)}$
\item[(ii)] eliminate $\alpha_{i+1}, \alpha_{i+2}, \dots, \alpha_{j-1}$
\item[(iii)] replace $\alpha_{\varphi(i)}=\alpha_j$ with $$\widehat{\alpha}_{(k,2)}\widehat{\alpha}_{k+1}\widehat{\alpha}_{k+2}\cdots \widehat{\alpha}_{j-1}\widehat{\widehat{\alpha}}_{i+1}\widehat{\widehat{\alpha}}_{i+2}\cdots \widehat{\widehat{\alpha}}_{k-1}\widehat{\widehat{\alpha}}_{(k,1)}\widehat{\alpha}^-_{(i,1)}$$
    where $\widehat{\alpha}_{(k,1)}=\widehat{\alpha}_k|_{[t^\ast_{k-1},b_1]}$, $\widehat{\alpha}_{(k,2)}=\widehat{\alpha}_k|_{[b_1,t^\ast_{k}]}$ (empty if $b_1=t^\ast_{k}$), and $\widehat{\widehat{\alpha}}_{\bullet}\equiv \widehat{\alpha}_{\bullet}$
\item[(iv)] replace (if defined) $\widehat{\alpha}_{\varphi(k)}, \widehat{\widehat{\alpha}}_{\varphi(k)}, \widehat{\alpha}_{\varphi^{-1}(k)}$, or $\widehat{\widehat{\alpha}}_{\varphi^{-1}(k)}$  with $\widehat{\alpha}^-_{(k,2)}\widehat{\alpha}^-_{(k,1)}$
\end{itemize}
(Note that there is no $\widehat{\alpha}_{(i,2)}$ and no $\widehat{\alpha}_{(j,2)}$.)
\vspace{10pt}

We then make the necessary adjustments to the semicircles, without introducing intersections: (See Figure~\ref{cycle}.)

\begin{itemize}
\item[(a)] We may assume that $S_i$ was already positioned so
as to serve as $\widehat{S}_{(i,1)}$,  matching $\alpha_{(i,1)}$ with $\widehat{\alpha}^-_{(i,1)}$.

\item[(b)] Take all semicircles over the domain of $\alpha_{i+1}\alpha_{i+2}\cdots\alpha_{j-1}$ and shift them to the corresponding positions of $\widehat{\alpha}_{i+1}\widehat{\alpha}_{i+2}\cdots\widehat{\alpha}_{j-1}$.

\item[(c)] If one of the semicircles touches the domain of $\widehat{\alpha}_k= \widehat{\alpha}_{(k,1)}\widehat{\alpha}_{(k,2)}$ (so that the pair $\widehat{\alpha}^-_{(k,2)}\widehat{\alpha}^-_{(k,1)}$ was created in step (iv)),  split it into two parallel semicircles, serving as $\widehat{S}_{(k,1)}$ and  $\widehat{S}_{(k,2)}$ (where the latter is not needed if $\widehat{\alpha}_{(k,2)}$ is empty).

\item[(d)] Circulate all semicircles under $S_i$ into their final corresponding positions.
\end{itemize}

\begin{figure}[h!]
\includegraphics{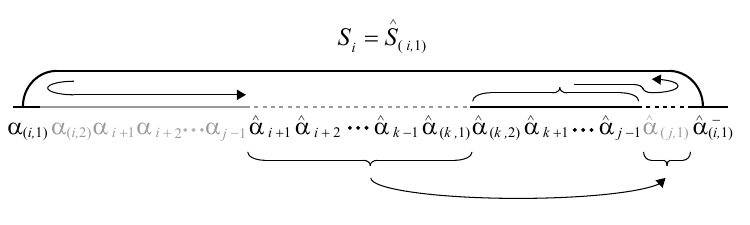}
\vspace{-20pt}

\caption{\label{cycle} Reposition, shift, split, and circulate semicircles}
\end{figure}
\pagebreak

{\bf Case~2.} Suppose $t^\ast_{j-1}<b_1<t^\ast_j$.
\vspace{10pt}

In this case,  it may take several rounds to transfer the subdivision points $\eta\cap [t_{j-1},b_1]$ to $[b_1,t_j]$.
Specifically, we proceed as follows (see Figure~\ref{Step3.2Case2}):

\begin{figure}[h!]
\hspace{-20pt} \includegraphics{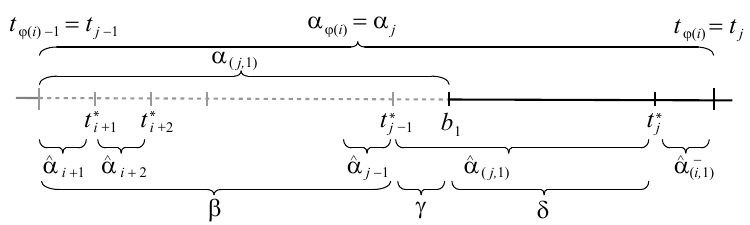}
\caption{\label{Step3.2Case2} Transfer of subdivision points (Step~3.2, Case2)}
\end{figure}

 Put $\beta=\alpha_j|_{[t_{j-1},t^\ast_{j-1}]}$,
 $\gamma=\alpha_j|_{[t^\ast_{j-1},b_1]}$ and $\delta=\alpha_j|_{[b_1,t^\ast_j]}$.
Note that $\gamma\delta=\widehat{\alpha}_{(j,1)}\equiv  \alpha_{(j,1)}=\beta\gamma$. Then $\beta\gamma\delta=\beta\beta_1\widehat{\gamma}$ with $\beta_1\equiv \beta$ and $\widehat{\gamma}\equiv \gamma$. If the terminal point of $\beta_1$ is less than or equal to $b_1$, then $\gamma=\beta_1\gamma_1$ for some $\gamma_1$, so that $\beta\gamma=\beta\beta_1\gamma_1$ and
 $\beta\gamma\delta=\beta\beta_1\beta_2\widehat{\gamma}_1$ with $\beta_2\equiv \beta_1$ and $\widehat{\gamma}_1\equiv \gamma_1$. If the terminal point of $\beta_2$ is still less than or equal to $b_1$, then $\gamma_1=\beta_2\gamma_2$ for some $\gamma_2$, so that $\beta\gamma=\beta\beta_1\beta_2\gamma_2$ and
 $\beta\gamma\delta=\beta\beta_1\beta_2\beta_3\widehat{\gamma}_2$ with $\beta_3\equiv \beta_2$ and $\widehat{\gamma}_2\equiv \gamma_2$. We continue inductively. Since $\alpha$ is continuous, there must be a first $\beta_s$ whose terminal point is greater than $b_1$. At that time,
$\beta\gamma=\beta\beta_1\beta_2\cdots\beta_{s-1}\gamma_{s-1}$ and
 $\beta\gamma\delta=\beta\beta_1\beta_2\cdots \beta_{s-1}\beta_{s}\widehat{\gamma}_{s-1}$ with $\beta_s\equiv\beta_{s-1}\equiv\cdots\equiv\beta_1\equiv \beta\equiv \widehat{\alpha}_{i+1}\widehat{\alpha}_{i+2}\cdots\widehat{\alpha}_{j-1}$ and $\widehat{\gamma}_{s-1}\equiv \gamma_{s-1}$.\linebreak
  Say, $\beta_s=\widehat{\widehat{\alpha}}_{i+1}\widehat{\widehat{\alpha}}_{i+2}\cdots\widehat{\widehat{\alpha}}_{j-1}
  =\widehat{\widehat{\alpha}}_{i+1}\widehat{\widehat{\alpha}}_{i+2}\cdots \widehat{\widehat{\alpha}}_{k-1} \widehat{\widehat{\alpha}}_{(k,1)}\widehat{\widehat{\alpha}}_{(k,2)}\widehat{\widehat{\alpha}}_{k+1}\cdots \widehat{\widehat{\alpha}}_{j-1}$, where we split $\widehat{\widehat{\alpha}}_{k}$ at $b_1$ in its domain. Since $\gamma_{s-1}\delta=\beta_s\widehat{\gamma}_{s-1}$, we have $$\delta\equiv\widehat{\widehat{\alpha}}_{(k,2)}\widehat{\widehat{\alpha}}_{k+1}\cdots \widehat{\widehat{\alpha}}_{j-1}\widehat{\widehat{\alpha}}_{i+1}\widehat{\widehat{\alpha}}_{i+2}\cdots \widehat{\widehat{\alpha}}_{k-1} \widehat{\widehat{\alpha}}_{(k,1)}$$ and we can proceed as in Case 1. (See Figure~\ref{Step3.2Return}.)

\begin{figure}[h!]
\hspace{.5in} \includegraphics{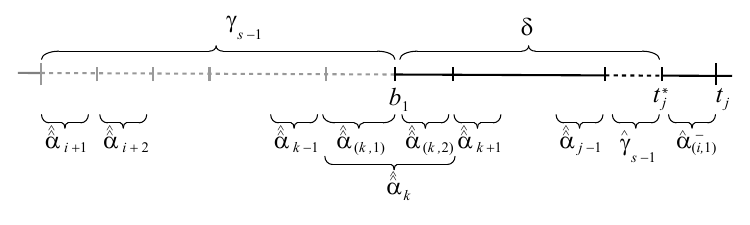}
\vspace{-25pt}

\caption{\label{Step3.2Return} Eventually, we can proceed as in  Case 1.}
\end{figure}

The stated goal was to eventually decrease $(r,n)$ in the lexicographical ordering. However,
the transformation described in Steps 1--3 may or may not decrease $r$ and potentially increases $n$ by as much as 2. The success analysis is now the same as in \cite{EF}. More specifically, if $\alpha_{(i,1)}$ is not empty, then $r$ decreases by 1. If both $\alpha_{(i,1)}$ and $\alpha_{(j,2)}$ are empty, then $r$ either decreases (by 1 or 2) or remains constant, but $n$ definitely decreases. So, assume now that $\alpha_{(i,1)}$ is empty and $\alpha_{(j,2)}$ is not empty. If $i=1$, then $r$ decreases from 1 to 0. So, assume now that $i>1$. If $\widetilde{\alpha}|_{[t_{i-2},t_{i-1}]\cup [b_1,b]}$ is a geodesic, then $r$ decreases by 1. Otherwise, $r$ remains constant, and $n$ does not increase. Verbatim as in the proof of \cite[Lemma~4.3]{EF} (using $\widetilde{\alpha}$ for $f$, $\alpha_i$ for $f_i$, $n$ for $s$, and $R$ for $F_o\cup F_e$) one shows that under repeated application of Steps 1--3, it is not possible for both $r$ and $n$ to remain constant indefinitely without an opportunity to eventually decrease $(r,n)$ in the lexicographical ordering.
\end{proof}

Recall that for a path-connected subset $U\subseteq X$, the subgroup $\pi(U)\leqslant \pi_1(X,x_0)$ is generated by all elements of the form $[\alpha_1\cdot \alpha_2\cdot \alpha_3]$ with $\alpha_1= \alpha_3^{-}$ and Im$(\alpha_2)\subseteq U$. If $X$ is one-dimensional, then each $\alpha_i$ can be chosen to be reduced.

\begin{corollary}\label{characterization}
Let $X$ be a path-connected one-dimensional separable metric space with basepoint $x_0\in X$. Suppose  $1\not=g\in \pi_1(X,x_0)$ and let $U\subseteq X$ be a path-connected subset. Then $g\in \pi(U)$ if and only if the reduced representative for $g$ has a cancellation relative to $U$.
\end{corollary}

\begin{proof}
The proof of the ``if'' part is  straightforward (using Remark~\ref{representation}, for example, and inducting on the number of vertices of the underlying tree structure). The ``only if'' part follows from Theorem~\ref{roundabout}.
\end{proof}

\begin{corollary}\label{strip}
For every $x\in B\subseteq \mathbb{D}$ and for every open neighborhood $V$ of $x$ in $\mathbb{D}$, there is a path-connected open neighborhood $U$ of $x$ in $\mathbb{D}$ such that $U\subseteq V$, $f(U)\subseteq U$, and with the following property: If $\beta$ is a  reduced loop in $\mathbb{D}$ based at $d_0$ with $[f\circ \beta]\in \pi(U)$, then $[\beta]\in \pi(U)$.

 In particular, the induced homomorphism $\overline{f}_\#$ in the following commutative diagram is injective:
 \[
 \xymatrix{\pi_1(\mathbb{D},d_0) \ar[r]^{f_\#} \ar[d] & \pi_1(\mathbb{D},d_0) \ar[d]\\
 \pi_1(\mathbb{D},d_0)/\pi(U) \ar[r]^{\overline{f}_\#} & \pi_1(\mathbb{D},d_0)/\pi(U)
 }
 \]
\end{corollary}

\begin{proof} Let $x=(t,0)\in B$. Let $n\in\mathbb{N}$ be given. Choose any $i\in \{0,1,\dots,2^{n}\}$ with $\frac{i-1}{2^{n}}<t<\frac{i+1}{2^{n}}$. Let $U$ be the path component of $x$ in $\mathbb{D}\cap \left((\frac{i-1}{2^{n}},\frac{i+1}{2^{n}})\times\mathbb{R}\right)$. Since $\mathbb{D}$ is locally path connected, $U$ is an open neighborhood of $x$ in $\mathbb{D}$. There are three basic cases for the shape of $U$: (i) $i=0$ or $i=2^n$; (ii) $0<i<2^n$ and $i$ is odd; (iii) $0<i<2^n$ and $i$ is even. In case (i), $f^{-1}(U)=U$. In cases (ii) and (iii), $f^{-1}(U)$ has two path-components, one of which is $U$ and the other is homeomorphic to an open interval. By choosing $n$ large enough, we can arrange for $U\subseteq V$.

Now, let $\beta$ be a reduced loop in $\mathbb{D}$ based at $d_0$ with $[f\circ \beta]\in \pi(U)$. We may assume that $\beta$ is not constant. Then $f\circ \beta$ is not constant. By Lemma~\ref{reduced},  $f\circ \beta$ is reduced. By Corollary~\ref{characterization}, there is a cancellation of $f\circ \beta$ relative to $U$. By Lemma~\ref{uniquelift}, the cancellation of $f\circ \beta$ relative to $U$ lifts to a cancellation of $\beta$ relative to $f^{-1}(U)$. If any one of the ``roundabout'' loops for $\beta$ (see Remark~\ref{representation}) lifts to a path component of $f^{-1}(U)$ which is homeomorphic to an open interval, then it can be further subdivided and incorporated into the cancellation pattern. Thus we arrive at a cancellation of $\beta$ relative to $U$. By Corollary~\ref{characterization}, $[\beta]\in \pi(U)$.
\end{proof}

\section{The mapping torus $\mathbb{M}$ of the shift-map $f$}

\begin{figure}[h!]
\includegraphics[scale=.6]{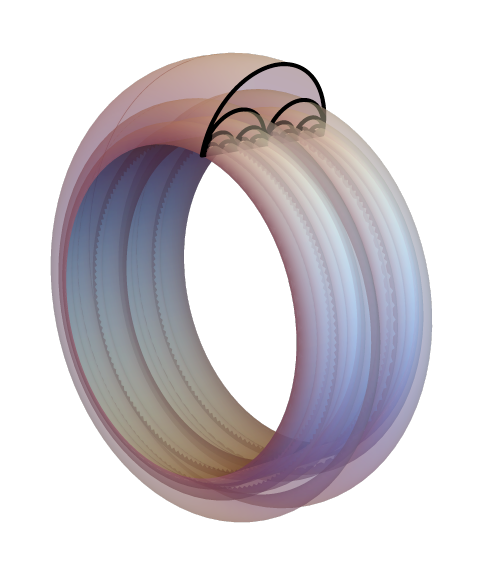}
\vspace{-10pt}

\caption{\label{M} The mapping torus $\mathbb{M}$ of  $f:\mathbb{D}\rightarrow \mathbb{D}$.}
\end{figure}

\noindent Let $\mathbb{M}$ be the mapping torus of the shift-map $f:\mathbb{D}\rightarrow \mathbb{D}$. That is, let $\mathbb{M}$ be the quotient space $\mathbb{D}\times [0,1]/\mathord{\sim}$ with $(x,1)\mathord{\sim} (f(x),0)$ for all $x\in \mathbb{D}$. Let $q:\mathbb{D}\times [0,1]\rightarrow \mathbb{M}$ denote the quotient map. Let $\xi:[0,1]\rightarrow \mathbb{M}$ be given by $\xi(t)=q(d_0,t)$. We identify $\mathbb{D}$ with $q(\mathbb{D}\times\{0\})\subseteq \mathbb{M}$. (See Figure~\ref{M}.)
\vspace{-10pt}

\begin{lemma}
$\mathbb{M}$ is a Peano continuum. Moreover, $\mathbb{M}$ is locally contractible at every $x\not\in q(B\times[0,1])$.
\end{lemma}

\begin{proof}
Note that $\mathbb{M}$ can be realized in $\mathbb{R}^3$ as the quotient of the Peano continuum $\mathbb{D}\times [0,1]$. The second statement is obvious.
\end{proof}

\begin{lemma}\label{inclusion}
Let $\iota:\mathbb{D}\hookrightarrow \mathbb{M}$ denote inclusion. Then the homomorphism $\iota_\#:\pi(\mathbb{D},d_0)\rightarrow \pi_1(\mathbb{M},d_0)$ is injective.
\end{lemma}

The proof of Lemma~\ref{inclusion} follows along the lines of the proof of \cite[Lemma 2.2]{BFi}. For later reference, we include the proof:

\begin{proof}[Proof of Lemma~\ref{inclusion}]
For each $i\in \mathbb{N}$, let $(X_i,x_i)$ be a copy of $(\mathbb{D},d_0)$ and let $f_i=f:X_i\rightarrow X_{i+1}$.
Let $\overline{\mathbb{M}}=\left(\sum_{i\in\mathbb{Z}} X_i\times[0,1]\right)/\mathord{\sim}$ denote the mapping telescope of the bi-infinite sequence \[\cdots  \stackrel{f_{-2}}{\longrightarrow} X_{-1}\stackrel{f_{-1}}{\longrightarrow} X_{0} \stackrel{f_{0}}{\longrightarrow} X_{1}\stackrel{f_{1}}{\longrightarrow} X_{2} \stackrel{f_{2}}{\longrightarrow} \cdots\]
where $(x,1)\mathord{\sim}(f_i(x),0)$ for all $i\in\mathbb{Z}$ and all $x\in X_i$.
We identify $X_i$ with the image of $X_i\times\{0\}$ in $\overline{\mathbb{M}}$.
Let $h:\overline{\mathbb{M}}\rightarrow \mathbb{M}$ be the canonical covering projection, mapping the image of $(x,t)\in X_i\times [0,1]$ in $\overline{\mathbb{M}}$ to the image of $(x,t)$ in $\mathbb{M}$.

Let $[\alpha]\in\pi_1(\mathbb{D},d_0)$ be such that $[\alpha]=1\in \pi_1(\mathbb{M},d_0)$. Then $\alpha$ lifts to a loop $\overline{\alpha}$ in  $X_0\subseteq \overline{\mathbb{M}}$ based at $x_0$ with $h\circ \overline{\alpha}=\alpha$ and $[\overline{\alpha}]=1\in\pi_1(\overline{\mathbb{M}},x_0)$.
Since any homotopy contracting $\overline{\alpha}$ in $\overline{\mathbb{M}}$ has compact image, and since for each $i\in\mathbb{Z}$, the image of $X_i\times[0,1]$ in $\overline{\mathbb{M}}$ canonically deformation retracts onto $f_i(X_i)\subseteq X_{i+1}$ along the second coordinate, there is a $k\in \mathbb{N}$ such that $f_k\circ f_{k-1}\circ\cdots\circ f_0\circ\overline{\alpha}$ contracts in $X_{k+1}$. However, $f_\#:\pi_1(\mathbb{D},d_0)\rightarrow \pi_1(\mathbb{D},d_0)$ is injective by Lemma~\ref{f-inj}.
Hence, $[\overline{\alpha}]=1\in \pi_1(X_0,x_0)$ and $[\alpha]=[h\circ \overline{\alpha}]=1\in \pi_1(\mathbb{D},d_0)$.
\end{proof}

\begin{lemma}\label{conjugate} For all loops $\gamma$ in $\mathbb{D}$ based at $d_0$, we have  $f_\#([\gamma])=[\xi]^{-1}[\gamma][\xi]\in \pi_1(\mathbb{M},d_0)$.
\end{lemma}

\begin{proof}
We use the notation from the proof of Lemma~\ref{inclusion}. The loop $\gamma$  lifts to a loop $\overline{\gamma}$ in $X_0$ based at $x_0$. Let $\overline{X_0\times [0,1]}$ be the image of $X_0\times [0,1]$ in $\overline{\mathbb{M}}$ and let $D:\overline{X_0\times [0,1]}\times[0,1]\rightarrow \overline{X_0\times [0,1]}$ be the canonical deformation retraction onto $f(X_0)\subseteq X_1\subseteq  \overline{\mathbb{M}}$.
Then $\gamma(t)=h\circ D(\overline{\gamma}(t),0)$, $f\circ \gamma(t)=h\circ D(\overline{\gamma}(t),1)$, and $\xi(t)=h\circ D(x_0,t)$. Hence, $[\gamma]=[\xi][f\circ\gamma][\xi]^{-1}\in \pi_1(\mathbb{M},d_0)$.
\end{proof}

\begin{theorem}\label{noninj} Let $\iota:\mathbb{D}\hookrightarrow \mathbb{M}$ denote inclusion. Then  $\iota_\#(\pi_1(\mathbb{D},d_0))$ is uncountable and $\iota_\#(\pi_1(\mathbb{D},d_0)) \leqslant \ker \Phi_{\mathbb{M},d_0}$. In particular, the homomorphism $\Phi_{\mathbb{M},d_0}:\pi_1(\mathbb{M},d_0)\rightarrow \check{\pi}_1(\mathbb{M},d_0)$ is not injective.
\end{theorem}

\begin{proof} Consider $\mathbb{D}\subseteq \mathbb{R}^2$ and $\mathbb{D}\subseteq \mathbb{M}\subseteq \mathbb{R}^3$. Let $d$ denote the standard metric for $\mathbb{R}^3$. Let $\gamma:([0,1],\{0,1\})\rightarrow (\mathbb{D},d_0)$ be a loop.
 Let $\left\{ {\mathcal U}_j, U^j_x,P_j,p_j, \psi_j, \varphi_j \right\}$ be a nice polyhedral expansion of $(\mathbb{M},d_0)$.

Say, $\gamma(t)=(x(t),y(t),0)$. Define $\alpha:[0,1]\rightarrow B$ by $\alpha(t)=(x(t),0,0)$. Since $B$ is contractible,  $\alpha$ is null-homotopic. Fix $j\in \mathbb{N}$. Choose a subdivision $0=t_0<t_1<\cdots<t_n=1$ of $[0,1]$ and elements $U_1, U_2, \dots, U_n\in {\mathcal U}_j$ such that $\alpha([t_{i-1}, t_i]) \subseteq  U_i$ for all $1 \leqslant i \leqslant n$. Choose $\epsilon>0$ so that the $\epsilon$-neighborhood of $\alpha([t_{i-1}, t_i])$ lies in $U_i$ for all $1 \leqslant i \leqslant n$, and the $\epsilon$-neighborhood of $\alpha_i(t_i)$ lies in $U_i\cap U_{i+1}$ for all $1 \leqslant i \leqslant n-1$. Choose $k\in\mathbb{N}$ such that $\sqrt{2}/2^{k+1}<\epsilon$.
Put $\delta=f^{k}\circ\gamma$.
Observe that for all $t\in [0,1]$, we have $d(\alpha(t),\delta(t))\leqslant 1/2^{k+1}$. Moreover, for all $1 \leqslant i \leqslant n-1$, there is a path $\kappa_i$ in $\mathbb{D}$ of diameter at most $\sqrt{2}/2^{k+1}$ from $\delta(t_i)$ to $\alpha(t_i)$. Put $\delta_i=\delta|_{[t_{i-1},t_i]}$ and $\alpha_i=\alpha|_{[t_{i-1},t_i]}$ for  $1 \leqslant i \leqslant n$. Put $\beta_1=\delta_1\kappa_1\alpha_1^-$, $\beta_i=\kappa_{i-1}^-\cdot \delta_i\cdot \kappa_i\cdot \alpha_i^-$ for $2\leqslant i\leqslant n-1$, and $\beta_n=\kappa_{n-1}\cdot \delta_n\cdot \alpha_n^-$.  Then $[\delta]=[\beta_1\cdot \alpha_1 \cdot \beta_2 \cdot \alpha_2\cdots \beta_n\cdot \alpha_n]$ with loops $\beta_i$ in $U_i$ and $[\alpha_1\cdot \alpha_2\cdots \alpha_n]=[\alpha]=1$. By Lemma~\ref{loops},  $[\delta]\in \pi(U_1)\pi(U_2)\cdots \pi(U_n)$. By Lemma~\ref{conjugate}, $[\gamma]=[\xi]^k[\delta][\xi]^{-k}\in \pi(U_1)\pi(U_2)\cdots \pi(U_n)$. By Proposition~\ref{ker}, $[\gamma]\in \ker\Phi_{\mathbb{M},d_0}$. Hence, $\iota_\#(\pi_1(\mathbb{D},d_0))\leqslant \ker \Phi_{\mathbb{M},d_0}$.

By Lemma~\ref{inclusion}, $\iota_\#(\pi_1(\mathbb{D},d_0))$ is isomorphic to $\pi_1(\mathbb{D},d_0)$, which is uncountable \cite{CC}.
\end{proof}

\begin{remark}\label{ncl} It follows from Theorem~\ref{noninj} that $h:\overline{\mathbb{M}}\rightarrow \mathbb{M}$ is a (categorical) universal covering projection and that $\ker \Phi_{\mathbb{M},d_0}$ equals $h_{\#}(\pi_1(\overline{\mathbb{M}},x_0))$, which in turn equals the normal closure of $\iota_\#(\pi_1(\mathbb{D},d_0))$ in $\pi_1(\mathbb{M},d_0)$. (Indeed, given $g\in \ker \Phi_{\mathbb{M},d_0}$ and any covering projection $r:(E,e_0)\rightarrow (\mathbb{M},x_0)$, we may choose $i\in \mathbb{N}$ such that each $U\in {\mathcal U}_i$ is evenly covered by $r$, so that $g\in r_{\#}(\pi_1(E,e_0))$ by Proposition~\ref{ker}. In turn, if $g\in h_{\#}(\pi_1(\overline{\mathbb{M}},x_0))$, then, as in the proof of Lemma~\ref{rewrite} below, $g=[\xi]^k[\gamma][\xi]^{-k}$ for some $k\in \mathbb{N}$ and some loop $\gamma$ in $\mathbb{D}$.)
\end{remark}

\begin{lemma}\label{rewrite} Let $g\in \pi_1(\mathbb{M},d_0)$ and $x\in B\subseteq \mathbb{D}$.
Suppose \[g=\prod_{i=1}^n [\beta_i\cdot \gamma_i\cdot \beta_i^-],\] where each $\beta_i$ is a path in $\mathbb{M}$ from $\beta_i(0)=d_0$ to $\beta_i(1)=x$, and each $\gamma_i$ is a loop in $\mathbb{D}$ at $x$.
Then for every $m\in\mathbb{N}$ there are $k>m$, $s_i\in\{0,1,\dots,2k\}$, and paths $\delta_i$ in $\mathbb{D}$ from $\delta_i(0)=d_0$ to $\delta_i(1)=x$ such that \[g =[\xi]^k\left(\prod_{i=1}^n[\delta_i\cdot (f^{s_i}\circ \gamma_i)\cdot \delta_i^-]\right) [\xi]^{-k}.\]
\end{lemma}

\begin{proof}
Again, we use the notation from the proof of Lemma~\ref{inclusion}. Let $\overline{\beta}_1$ be the lift of $\beta_1$ with $\beta_1=h\circ \overline{\beta}_1$ and $\overline{\beta}_1(0)=x_0$. Then $\overline{\beta}_1(1)=y$ for some $y\in X_j$ and $j\in \mathbb{Z}$ such that $h(y)=x$.  Let $\overline{\gamma}_1$ be the lift of $\gamma_1$ with $\gamma_1=h\circ \overline{\gamma}_1$ and $\overline{\gamma}_1(0)=y$. Then $\overline{\gamma}_1$ is a loop in $X_j$ based at $y$ and $\overline{\beta_1}\cdot \overline{\gamma_1}\cdot (\overline{\beta_1})^-$ is a loop in $\overline{\mathbb{M}}$ based at $x_0$. Since the image of $\overline{\beta_1}\cdot \overline{\gamma_1}\cdot (\overline{\beta_1})^-$ is compact, there is a $K_1\in \mathbb{N}$ such that for every $k>K_1$, the loop $\overline{\beta_1}\cdot \overline{\gamma_1}\cdot (\overline{\beta_1})^-$ is contained in the subset $\overline{\mathbb{M}}_k$ of $\overline{\mathbb{M}}$ corresponding to the mapping telescope of the finite subsequence \[X_{-k}\stackrel{f_{-k}}{\longrightarrow} X_{-k+1} \stackrel{f_{-k+1}}{\longrightarrow}  \cdots \stackrel{f_{-1}}{\longrightarrow} X_0 \stackrel{f_{0}}{\longrightarrow}\cdots \stackrel{f_{k-2}}{\longrightarrow} X_{k-1} \stackrel{f_{k-1}}{\longrightarrow} X_k.\]
In particular, $-k\leqslant j \leqslant k$. As in the proofs of the previous two lemmas, we may now apply successive deformation retractions of the image of $X_s\times [0,1]$ in $\overline{\mathbb{M}}$ onto $f(X_s)\subseteq X_{s+1}$  for $s=-k, -k+1, \dots, k-1$ and projecting them into $\mathbb{M}$ using $h$, to see that \[[\beta_1\cdot \gamma_1\cdot \beta_1^-]=[\xi]^k[\delta_1\cdot (f^{s_1}\circ \gamma_1)\cdot \delta_1^-][\xi]^{-k},\] where $s_1=k-j$ and $\delta_1$ is a path in $\mathbb{D}$ from $d_0$ to $x$. We similarly find $K_2, K_3, \dots, K_n\in \mathbb{N}$ such that for all $i\in \{1,2,\cdots,n\}$ and all $k>K_i$, there are $s_i\in \{0,1,\dots,2k\}$ and paths $\delta_i$ in $\mathbb{D}$ from $d_0$ to $x$ such that $[\beta_i\cdot \gamma_i\cdot \beta_i]=[\xi]^k[\delta_i\cdot (f^{s_i}\circ \gamma_i)\cdot \delta_i^-][\xi]^{-k}$.
To prove the lemma, choose $k>\max\{m,K_1,K_2,\dots,K_n\}$.
\end{proof}

\begin{theorem}
For every $x\in \mathbb{M}$ and  every $1\not=g\in\pi_1(\mathbb{M},d_0)$,  there exists a path-connected open neighborhood $W$ of $x$ in $\mathbb{M}$ such that \linebreak $g\not\in \pi(W)\leqslant \pi_1(\mathbb{M},d_0)$.
\end{theorem}

\begin{proof} Let $x\in \mathbb{M}$ and  $1\not=g\in\pi_1(\mathbb{M},d_0)$.
Suppose, to the contrary, that $g\in \pi(W)\leqslant \pi_1(\mathbb{M},d_0)$ for every path-connected open neighborhood $W$ of $x$ in $\mathbb{M}$. Then $x\in q(B\times [0,1])$. We may assume that $x\in B$. (Otherwise, the proof is similar.)
Given any open subset $V$ of $\mathbb{D}$, we define \[V^+=q((f^{-1}(V)\times (\frac{2}{3},1])\cup (V\times [0,\frac{1}{3}))).\] Then $V^+$ is an open subset of $\mathbb{M}$ which deformation retracts onto $V$.

Since $g\in \pi(\mathbb{D}^+)\leqslant \pi_1(\mathbb{M},d_0)$ and since $\mathbb{D}^+$ deformation retracts onto $\mathbb{D}$, it follows from Lemma~\ref{rewrite} that there is an $m\geqslant 0$ and a nontrivial loop $\beta$ in $\mathbb{D}$ based at $d_0$ such that $g=[\xi]^m[\beta][\xi]^{-m}\in \pi_1(\mathbb{M},d_0)$. We may assume that $\beta$ is reduced.

Since $\Phi_{\mathbb{D},d_0}:\pi_1(\mathbb{D},d_0)\rightarrow \check{\pi}_1(\mathbb{D},d_0)$ is injective, by Proposition~\ref{ker} and Remark~\ref{local}, there is a path-connected open neighborhood $U$ of $x$ in $\mathbb{D}$ such that $[\beta]\not\in \pi(U)\leqslant \pi_1(\mathbb{D},d_0)$. We may assume that $U$ has the properties stated in Corollary~\ref{strip}: $f(U)\subseteq U$ and whenever $\beta'$ is a  reduced loop in $\mathbb{D}$ based at $d_0$ with $[f\circ \beta']\in \pi(U)\leqslant\pi_1(\mathbb{D},d_0)$, then $[\beta']\in \pi(U)\leqslant \pi_1(\mathbb{D},d_0)$.

Consequently,  noting Lemma~\ref{reduced},
\begin{equation}\label{U}\tag{$\ast$}
f^s_\#([\beta])\not \in \pi(U)\leqslant \pi_1(\mathbb{D},d_0) \mbox{ for all } s\geqslant 0.
\end{equation}

Since $g\in \pi(U^+)\leqslant \pi_1(\mathbb{M},d_0)$, we have \[g=\prod_{i=1}^n [\beta_i\cdot \gamma_i\cdot \beta_i^-]\] where
each $\beta_i$ is a path in $\mathbb{M}$ from $\beta_i(0)=d_0$ to $\beta_i(1)=x$, and each $\gamma_i$ is a loop in $U^+$ at $x$. We may assume that each $\gamma_i$ lies in $U$.
By Lemma~\ref{rewrite}, we have
\[g =[\xi]^k\left(\prod_{i=1}^n[\delta_i\cdot\gamma_i'\cdot \delta_i^-]\right) [\xi]^{-k}\]
with $k>m$, $s_i\in\{0,1,\dots,2k\}$, paths $\delta_i$ in $\mathbb{D}$ from $\delta_i(0)=d_0$ to $\delta_i(1)=x$ and loops $\gamma_i'=f^{s_i}\circ \gamma_i$ in $U$ based at $x$.

Then
\[[\xi]^m[\beta][\xi]^{-m}=[\xi]^k\left(\prod_{i=1}^n[\delta_i\cdot\gamma_i'\cdot \delta_i^-]\right) [\xi]^{-k},\]
which by Lemma~\ref{conjugate} yields
\[f_\#^{k-m}([\beta])=\prod_{i=1}^n[\delta_i\cdot\gamma_i'\cdot \delta_i^-]\in \pi(U)\leqslant \pi_1(\mathbb{M},d_0).\]
Then, by Lemma~\ref{inclusion},
\[f_\#^{k-m}([\beta])=\prod_{i=1}^n[\delta_i\cdot\gamma_i'\cdot \delta_i^-]\in \pi(U)\leqslant \pi_1(\mathbb{D},d_0),\]
which contradicts (\ref{U}).\end{proof}

\begin{corollary}\label{main}
Let $\left\{ {\mathcal U}_i, U^i_x,P_i,p_i, \psi_i, \varphi_i \right\}$ be a nice polyhedral expansion of $(\mathbb{M},d_0)$. Then $\bigcap_{i\in \mathbb{N}} \pi(U^i_x)=1$  for all $x\in \mathbb{M}$.
\end{corollary}

\begin{theorem}\label{T1} $\pi_1(\mathbb{M},d_0)$ is T$_1$ (in the quotient of the compact-open topology).
\end{theorem}

\begin{proof} Fix a metric $d$ for $\mathbb{M}$. Let $1\not=[\alpha]\in \pi_1(\mathbb{M},d_0)$. We wish to find an $\epsilon>0$ such that there is no null-homotopic loop  $\alpha':([0,1],\{0,1\})\rightarrow (\mathbb{M},d_0)$ with $d(\alpha(t),\alpha'(t))<\epsilon$ for all $t\in [0,1]$.

We use the notation from the proof of Lemma~\ref{inclusion}.
Recall that the covering projection $h:\overline{\mathbb{M}}\rightarrow \mathbb{M}$ has continuous lifting of paths \cite[2.7.8]{S}.
Consider the space $\Omega(\mathbb{M},d_0)$  of all loops in $\mathbb{M}$ based at $d_0$, in the compact-open topology.
For $\alpha'\in \Omega(\mathbb{M},d_0)$, let $\overline{\alpha}':[0,1]\rightarrow \overline{\mathbb{M}}$ denote the lift of $\alpha'$ with $h\circ \overline{\alpha}'=\alpha'$ and $\overline{\alpha}'(0)=x_0$.
For $\epsilon>0$, define $N_\epsilon(\alpha)=\{\alpha'\in \Omega(\mathbb{M},d_0) \mid d(\alpha(t),\alpha'(t))<\epsilon$ for all $t\in [0,1]\}$.
Since the fiber $h^{-1}(d_0)$ is discrete, there is an $\epsilon_1>0$  such that for every $\alpha'\in N_{\epsilon_1}(\alpha)$ we have $\overline{\alpha}'(1)=\overline{\alpha}(1)$.
In particular, if $\overline{\alpha}(1)\not=x_0$, then for every $\alpha'\in N_{\epsilon_1}(\alpha)$, we have $1\not=[\alpha']\in\pi_1(\mathbb{M},d_0)$. So, let us now assume that $\overline{\alpha}$ is a loop.

For $k\in \mathbb{N}$, let $\overline{\mathbb{M}}_k\subseteq \overline{\mathbb{M}}$ be as in the proof of Lemma~\ref{rewrite} and let $D_k:\overline{\mathbb{M}}_k\rightarrow X_k$ denote the end of the canonical deformation retraction along the second coordinate.
Choose $k\in \mathbb{N}$ and $\epsilon_2\in(0,\epsilon_1)$ such that for every $\alpha'\in N_{\epsilon_2}(\alpha)$, we have $\overline{\alpha}'([0,1])\subseteq \overline{\mathbb{M}}_k$.
Then for each $\alpha'\in N_{\epsilon_2}$, we have $[\overline{\alpha}']=[\zeta\cdot (D_k\circ \overline{\alpha}') \cdot \zeta^-]\in\pi_1(\overline{\mathbb{M}},x_0)$ with $h\circ \zeta=\xi\cdot\xi\cdots \xi$, a $k$-fold concatenation.
Since $1\not=[\alpha]\in \pi_1(\mathbb{M},d_0)$, we have $1\not=[\overline{\alpha}]\in \pi_1(\overline{\mathbb{M}},x_0)$, so that  $1\not=[D_k\circ \overline{\alpha}]\in \pi_1(X_k,x_k)$.
 Since $(X_k,x_k)$ is a copy of $(\mathbb{D},d_0)$, it follows from the proof of Lemma~\ref{f-inj} that there is an $n\in\mathbb{N}$ with $1\not=r_{n\#}([D_k\circ \overline{\alpha}])\in \pi_1(E_n,d_0)$. Since the function $N_{\epsilon_2}(\alpha)\rightarrow \Omega(E_n,d_0)$ given by $\alpha'\mapsto r_n\circ D_k \circ \overline{\alpha}'$ is continuous and since $E_n$ is a finite graph, there is an $\epsilon_3\in (0,\epsilon_2)$ such that for every $\alpha'\in N_{\epsilon_3}(\alpha)$ we have $r_{n\#}([ D_k\circ\overline{\alpha}'])=r_{n\#}([D_k\circ\overline{\alpha}])\not=1$, so that $1\not=[D_k\circ \overline{\alpha}']\in \pi_1(X_k,x_k)$, $1\not=[\overline{\alpha}']\in\pi_1(\overline{\mathbb{M}},x_0)$  and $1\not=[\alpha']=h_\#([\overline{\alpha}'])\in \pi_1(\mathbb{M},d_0)$.
\end{proof}

\section{Comparison with other local properties}\label{compare}

\noindent Finally, we want to insert our results into the context of other local properties of fundamental groups studied elsewhere, for various reasons.

First, let us generalize the definition of the subgroup $\pi(U)\leqslant \pi_1(X,x_0)$ from Section~\ref{HPol} by not requiring that $U$ be path connected.

\begin{definition}
Let $X$ be a path-connected space with basepoint $x_0$. Let $Cov(X)$ denote the collection of all open covers of $X$. For $x\in X$, let ${\mathcal T}_x$ denote the collection of all open neighborhoods of $x$ in $X$.
 The subgroup \[\pi^s(X,x_0)=\bigcap_{{\mathcal U}\in Cov(X)} \left< \pi(U)\mid U\in {\mathcal U}\right>\leqslant \pi_1(X,x_0)\] is  referred to as the {\em Spanier group} of $(X,x_0)$ \cite{S}.
 We call the subgroup \[\pi^x(X,x_0)=\bigcap_{U\in {\mathcal T}_x} \pi(U)\leqslant \pi^s(X,x_0)\] the {\em point group} of $(X,x_0)$ at $x$.
\end{definition}

\begin{definition}
Let $X$ be a path-connected space with basepoint $x_0$. We call $X$ {\em completely homotopically Hausdorff} (CHH) if $\pi^x(X,x_0)=1$ for every $x\in X$.
\end{definition}

Observe that CHH is stronger than the property of being {\em strongly homotopically Hausdorff} (SHH) from \cite{CMRZZ}, which stipulates that for every  $x\in X$ and for every essential loop $\gamma$ in $X$, there be a $U\in {\mathcal T}_x$ such that $\gamma$ cannot be freely homotoped into $U$.

In turn, the property of being {\em homotopically Hausdorff} (HH) from \cite{CC} only requires that for every $x\in X$ and for every $1\not=g\in \pi_1(X,x)$, there is a $U\in {\mathcal T}_x$ such that $g\not=[\beta]$ for every loop $\beta$ in $U$ based at $x$.

The property of being homotopically path Hausdorff from \cite{FRVZ}, on the other hand, is a path-local  strengthening of HH, which is equivalent to $\pi_1(X,x_0)$ being T$_1$ for locally path-connected $X$ \cite{BFa2015}:

\begin{definition} A path-connected space $X$ is called
{\em homotopically path Hausdorff} (HPH) if for every two paths $\alpha, \beta : [0,1] \rightarrow X$ with $\alpha(0) = \beta(0)$ and
$\alpha(1) = \beta(1)$ such that $\alpha \cdot \beta^-$ is not null-homotopic, there is a partition $0 = t_0 < t_1 < \cdots <
t_n = 1$ of $[0,1]$ and open subsets $U_1, U_2, \dots , U_n$ of $X$ with $\alpha([t_{i-1}, t_i]) \subseteq  U_i$ for all $1 \leqslant i \leqslant n$ and with the property that if $\gamma : [0, 1] \rightarrow X$ is any path with $\gamma([t_{i-1}, t_i]) \subseteq  U_i$ for all $1 \leqslant i \leqslant n$ and with $\gamma(t_i) = \alpha(t_i)$ for all $0 \leqslant i \leqslant n$, then $\gamma \cdot \beta^-$ is not null-homotopic.
\end{definition}

The space $\mathbb{M}$ then finds its proper place in the following diagram of implications (top row) and non-implications (bottom row), where $\mathbb{P}$ is the space from \cite{BFi2020} and $\mathbb{T}$ is the space from \cite{BFi}:
\vspace{10pt}

\[
\xymatrix{ \text{HPH} \ar@/^-1.5pc/[r]_{\mathbb{M}}|{\SelectTips{cm}{}\object@{x}}  & \pi^s=1  \ar@/_1.5pc/[l] \ar@/^1.5pc/[r] & \text{CHH} \ar@/^1.5pc/[r] \ar@/_-1.5pc/[l]^{\mathbb{M}}|{\SelectTips{cm}{}\object@{x}}    & \text{SHH} \ar@/^1.5pc/[r] \ar@/_-1.5pc/[l]^{\mathbb{P}}|{\SelectTips{cm}{}\object@{x}}  & \text{HH} \ar@/_-1.5pc/[l]^{\mathbb{T}}|{\SelectTips{cm}{}\object@{x}}  }
\]

For further local properties of fundamental groups and their relationships, see \cite{BFi2020a}.

\end{document}